\documentclass[12pt]{amsart}
\usepackage{amssymb}
\usepackage{amsfonts}
\usepackage{amscd}
\usepackage[T1]{fontenc}

\usepackage{xy}
\input xy
\xyoption{all}
\usepackage{pdflscape}
\usepackage{hyperref}

\usepackage[english]{babel}

\def\A{\mathbb{A}}

\def\R{\mathbb{R}}
\def\C{\mathbb{C}}
\def\F{\mathbb{F}}

\def\GL{\mathrm{GL}}

\def\cc{\mathrm{cc}}

\let\myacute=\'

\def\ra{\longrightarrow}

\def\<{\langle}
\def\>{\rangle}
\def\N{\mathbb{N}}
\def\Z{\mathbb{Z}}

\def\cL{\mathcal{L}}

\def \begindm {\begin{displaymath}}

\def \enddm {\end{displaymath}}

\def \suchthat {\ensuremath{:}}

\def\C{\mathbb{C}}
\def\R{\mathbb{R}}
\def\Q{\mathbb{Q}}
\def\F{\mathbb{F}}

\def\cL{\mathcal{L}}

\def\cM{\mathcal M}
\def\cO{\mathcal{O}}

\newtheorem{thm}{Theorem}[section]
\newtheorem{lemma}{Lemma}[section]
\newtheorem{cor}{Corollary}[section]

\newtheorem{ex}{Example}[section]

\numberwithin{equation}{section}
\newtheorem{prob}{Problem}[section]
\newtheorem{note}{Note}[section]

\long\def\symbolfootnote[#1]#2{\begingroup\def\thefootnote{\fnsymbol{footnote}}\footnote[#1]{#2}\endgroup}

\title[$p$-adic integrals, Euler products, and zeta functions]{$p$-adic model theory, $p$-adic integrals, Euler products, and zeta functions of groups}
\author[J. Derakhshan]{Jamshid Derakhshan}
\address{St Hilda's College, University of Oxford, Cowley Place, Oxford OX4 1DY, UK}
\email{derakhsh@maths.ox.ac.uk}

\begin{document}



\begin{abstract} 
We give a survey of Denef's rationality theorem on $p$-adic integrals, its uniform in $p$ versions, the relevant model theory, and a number of applications to counting subgroups of finitely generated nilpotent groups and conjugacy classes in congruence quotients of Chevalley groups over rings of integers of local fields. We then state results on analytic properties of 
Euler products of such $p$-adic integrals over all $p$, and an application to counting conjugacy classes in congruence quotients of certain algebraic groups over the rationals. 
We then briefly discuss zeta functions arising from definable equivalence relations and $p$-adic elimination of imginaries, which have applications to counting representations of groups.

\end{abstract}

\maketitle

\tableofcontents

\section{\bf Introduction}\label{sec-introduction}
In 1984, Denef \cite{Denefrationality} 
proved a remarkable theorem on rationality of certain $p$-adic integrals. The integrals were analogous to the local zeta functions of Igusa \cite{igusa-book}, but their domain of integration were definable sets in the sense of logic. While 
Igusa had used the rationality of his local zeta functions to prove a conjecture of Borevich and Shafarevich on rationality of 
a Poincare series counting points on a variety modulo powers of $p$, Denef used his theorem to prove rationality of a 
Poincare series counting points modulo powers of $p$ that lift to a $p$-adic point. This proved a conjecture of Serre \cite{serre-conj}. 

Denef's theorem relies heavily on a quantifier elimination theorem of Macintyre for the $p$-adic numbers in the Macintyre language \cite{Macintyre1}. He combined this with Hironaka's embedded resolution of singularities for algebraic varieties and a model-theoretic cell decomposition theorem extending a result of Paul Cohen that gave a decision procedure for $p$-adic fields \cite{Cohen}. 

Several application followed, to number theory and group theory. Here we focus on the latter (for the former see \cite{adeles-surv}). Grunewald, Segal, and Smith \cite{GSS} introduced a zeta function counting subgroups of finite index in a finitely generated nilpotent group, and proved that it admits a factorization as an Euler product of $p$-adic integrals that are in turn zeta functions counting $p$th power index subgroups, and are rational by Denef's theorem. The subject of subgroup growth 
has been much developed in this connection (see the book by Lubotzky and Segal \cite{alex-dan-book} and Section \ref{growth}).

du Sautoy introduced a zeta function counting conjugacy classes in congruence quotients in in $GL_n(\Z_p)$ 
and used Denef's theorem to prove its rationality.

In joint work with Mark Berman, Uri Onn, and Pirita Paajanen \cite{BDOP} we defined the conjugacy class zeta function for all Chevalley groups, and proved that it only depends on the residue field provided the residue characteristic is large enough. 
For example $GL_n(\Z_p)$ and $GL_n(\F_p[[t]])$ have the same conjugacy class zeta function. This implies that 
for any $k\geq 1$, the groups $GL_n(\Z_p/p^k\Z_p)$ and $GL_n(\F_p[[t]]/t^k)$ have the same number of conjugacy classes, once $p$ is larger than some constant (depending only on $GL_n$). Same is true for any Chevalley group. 

This is an example of a model-theoretic transfer principle that was originally proved for truth of a sentence by Ax-Kochen, and later for an identity of motivic or definable integral by Denef, Loeser, and Cluckers in various forms (for example in \cite{DL}, \cite{CL2}). A version for integrals is also proved in \cite{BDOP} directly without using motivic integration. These are discussed in Section \ref{conj}.

In \cite{Denefrationality}, to apply the results on $p$-adic integration to prove Serre's conjecture, 
Denef proved "$p$-adic definability" of the Serre series. Analogously, in \cite{BDOP} we need to prove $p$-adic 
definability of the local conjugacy class zeta function. These proofs are presented in Sections \ref{serre} and \ref{def-conj} respectively. They provide the connection between model theoretic $p$-adic integrals with the generating functions in algebra and number theory, and open the door to apply $p$-adic model theory and $p$-adic integration the problems. We hope that the methods in these proofs can be used in various other generating functions and Poincare series. 

In \cite{zeta1}, Euler products of the above mentioned definable $p$-adic integrals are appropriately defined and it is proved that they admit meromorphic continuation beyond their abscissa of convergence and that one can get information on the poles of the Euler product. This is presented in Section \ref{euler} and the relevant model theory in Section \ref{qe}.

In \cite{zeta1}, this result combined with the uniformities of the local conjugacy class zeta functions are applied to a global conjugacy class zeta function for Chevalley groups with strong approximation (e.g. $SL_n$ for 
$n\geq 2$ and $Sp_{2g}$ for $g\geq 1$) to count conjugacy classes in congruence quotients over a number field, e.g. to get asymptotic formulas for the numbers of conjugacy classes in $SL_n(\Z/m\Z)$ as $m\rightarrow \infty$. This is discussed in Section \ref{global-conj}.

In the final section we discuss the subject of $p$-adic imaginaries (i.e. equivalence classes of a definable equivalence relation) and representation growth of a group that studies the 
numbers of representations of the group in any given degree using the representation zeta function. 

Hrushovski-Martin-Rideau \cite{HMR} proved that the $\Q_p$ have an elimination of imaginaries that is uniform in $p$ in 
the geometric language which has sorts for the spaces of lattices $GL_n(\Q_p)/GL_n(\Z_p)$ for all $n$. They used this 
to prove an extension of Denef's rationality theorem to Poincare series counting classes of a parametrized family of 
definable equivalence relations. They applied this to prove rationality of zeta functions counting (iso-twist classes of) representations of finitely generated nipotent groups. 

Moreover, Avni \cite{Nir} has proved results on representation growth for arithmetic lattices. These are lattices in semi-simple Lie groups of higher rank. If the lattice satisfies the congruence subgroup property, then the number of representations in each degree is finite, and one can define the representation zeta function. Avni proves results on the growth of the number of representations in each degree. He proves the representation zeta function has an Euler product factorization, over primes $p$, of certain 
$p$-adic integrals which are over definable sets and are a generalization of the integrals of Denef, and proves they are rational functions.

\section{\bf Acknowledgments}

I am very grateful to Nir Avni, Mark Berman, Raf Cluckers, Ehud Hrushovski, Francois Loeser, Marcus du Sautoy, Dan Segal, Angus Macintyre, Uri Onn, and Boris Zilber for many helpful discussions.

\section{\bf $p$-adic numbers and measures}\label{ssec-euler}

Let $p$ be a prime number.  Let $v_p(x)$ denote the 
$p$-adic valuation on $\Q$ defined by $v_p(x)=k$ if $x=p^k a/b$ where $ab$ is not divisible by $p$. This gives a $p$-adic 
absolute value that is non-Archimedean (or ultrametric) defined by $|x|=p^{-v(x)}$. The completion of $\Q$ under 
this absolute value is the field of $p$-adic numbers $\Q_p$. $v_p(x)$ admits a unique extension to a valuation of $\Q_p$, that we still denote by $v_p(x)$. Any element $a\in \Q_p$ can be written in the form 
$\sum_{i\geq -k} c_i p^i$, where $k$ is a non-negative integer and $c_i$ are either in $\{0,\dots,p-1\}$ or from the 
Teichmuller group. See \cite{cassels} for details. 

Let $K$ be a finite extension of $\Q_p$.
The $p$-adic valuation $v_p(x)$ admits a unique extension to $K$, that we denote by $v(x)$. 
Let $\cO_K=\{x\in K: v(x)\geq 0\}$ denote the valuation ring of $K$. 
Let $\pi$ be a uniformizing element of $K$, i.e. element of $K$ of least positive value, then 
any element $x\in K$ can be written as $x=\pi^k u$, where $k\in \Z$ and $u$ is a unit of $\cO_K$. $\cO_K$ is a local ring with 
unique maximal ideal $\cM_K=\{x\in K: v(x)>0\}$ generated by a uniformizing element. The field $k=\cO_K/\cM_K$ is called the 
residue field of $K$, and is finite of cardinality $q$ which is a power of $p$. $K$ carries an absolute value defined by 
$|x|=q^{-v(x)}$. See \cite{cassels}, \cite{ramak}.

The field $K$, hence the additive group $(K^n,+)$, is locally compact, thus carries an $\R$-valued Borel measure $\mu_n$ that is unique up to multiplication by a constant and invariant under translation. Since $\cO_K^n$ is compact, we normalize $\mu_n$ such that $\cO_K^n$ has volume $1$. We shall also denote this normalized additive measure on $K^n$ by 
$dx=dx_1\dots dx_n$.

Note that $$\mu_n(a+\pi^m \cO_K^n)=q^{-mn},$$
and for any measurable set $A$ and $\lambda \in K$, 
$$\mu_n(\lambda A)=|\lambda|\mu_n(A).$$ 
More generally, for any $g\in GK_n(K)$,
$$\mu_n(gA)=|\mathrm{det}(g)|\mu_n(A).$$

$K$ is a locally compact non-Archimedean field with finite residue field of characteristic zero. The topology is given by the metric $d(x,y)=|x-y|$. Conversely every such field is a finite extension of $\Q_p$ for some $p$. In characteristic $p>0$, such fields are exactly fields of Laurent series $\F_q((t))$. See \cite{ramak}

\section{\bf $p$-adic integration on analytic manifolds}

For the basic theory of $K$-manifolds and $K$-analytic functions see Serre's book \cite{serre-LGLA}.  If $f$ is a $K$-analytic function on a measurable set $A\subseteq K^n$, one defines
$$\int_{A}|f(x)|dx=\sum_{m\in \Z} \mu_n(\{x\in K^n: v(f(x))=m\})q^{-m},$$
assuming that it is convergent in $\R$, and 
$$\int_{A}|f(x)|^s dx=\sum_{m\in \Z} \mu_n(\{x\in K^n: v(f(x))=m\})q^{-ms},$$
when it is convergent, where $s$ is a complex variable.

\begin{ex} Let $n=1$. Let $s\in \R_{>0}$. Then
$$\int_{\{x\in \cO_K: v(x)\geq m\}} |x|^s=\sum_{j\geq m} q^{-sj} \int_{\{x\in \cO_K: v(x)=j\}} dx=\sum_{j\geq m} q^{-sj}(g^{-j}-q^{-j-1})$$
$$=(1-q^{-1})q^{-(s+1)m}/(1-q^{-(s+1)}).$$
\end{ex}
A fundamental fact is the following $p$-adic change of variables formula. 
\begin{thm}\cite{weil-adeles-gps},\cite{igusa-book} Let $U$ be an open subset of $K^n$ and $f_1,\dots,f_n$ $K$-analytic functions on $U$. Suppose that $f=(f_1,\dots,f_n): U\rightarrow K^n$ is a $K$-analytic isomorphism between $U$ and an open subset 
$V\subseteq K^n$. Then for every integrable function $\varphi$ on $V$,
$$\int_V \varphi \hspace{0.1cm} \mu_n|_{V}=\int_U (\varphi \circ f) |\Omega(f_1,\dots,f_n)| \hspace{0.1cm} \mu_n|_{U},$$
where $\Omega(f_1,\dots,f_n)$ is the determinant of the Jacobian matrix of $f$.\end{thm}
If $X$ is an $n$-dimensional smooth $K$-analytic manifold, then any $K$-analytic $n$-differential form $\omega$ on $X$ gives rise to a measure $\mu_{\omega}$ on $X$ defined as follows. Let $\{(U,\phi_U\}$ be an atlas of $X$. Let 
$(\phi_U^{-1})^*\omega|_{U}=f_U(x) dx$. If $A$ is included in some $U$, then one defines
$$\mu_{\omega}(A)=\int_{\phi_U(A)}|f_U(x)|dx.$$
By the $p$-adic change of variables formula, this can be extended to all measurable $A$ independent of choice of 
an atlas.

\section{\bf Conjectures of Borevich-Shafarevich and Serre}
Let $f_1(x),\dots,f_r(x)$ be polynomials in $m$ variables $x=(x_1,\dots,x_n)$ over $\Z_p$.  For $n\in \N$, let $N_n$ denote the 
number of elements in the set
$$\{x \in (\Z_p/p^n \Z_p)^m: f_i(x)\equiv 0 \ \mathrm{mod} \ p^m, \mathrm{for} \ i=1,\dots,r\}$$
and $\tilde{N}_n$ the number of elements in the set
$$\{x\in (\Z_p/p^n \Z_p)^m: \exists y \in \Z_p^m, y\equiv x \ \mathrm{mod} \ p^n, f_i(x)=0, \mathrm{for} \ i=1,\dots,r\}.$$
Consider the following Poincare series

$$\tilde{P}(T)=\sum_{n\geq 0} \tilde{N}_n T^n, \ P(T)=\sum_{n\geq 0} N_n T^n.$$
Borevich and Shafarevich \cite[page 6]{bor-shaf-book} conjectured that $P(T)$ is a rational function of $T$. This was proved by 
Igusa \cite{igusa-cx-pow1} when $r=1$ using resolution of singularities (see also Igusa's book \cite{igusa-book}), and by Meuser \cite{Meuser} for any $r$ using Igusa's method. 

Serre [Section 3]\cite{serre-conj} and Oesterle \cite{oesterle} 
conjectured that $\tilde{P}(T)$ is a rational function of $T$. This was proved 
by Denef \cite{Denefrationality}. 
\begin{thm}[Denef \cite{Denefrationality}] \label{denef-thm1} The series $\tilde{P}(T)$ is a rational function of $T$.
\end{thm}
Denef gave two proofs of this theorem, one using a theorem of Hironaka on resolution of singularities of an algebraic variety (see \cite{igusa-book}), and another using instead and a cell decomposition theorem \cite{Denefrationality} which extends work of Paul Cohen \cite{Cohen}. In both proofs, a theorem of Macintyre \cite{Macintyre1} giving quantifier elimination theorem for $\Q_p$ in is crucially used. Denef also gave a new proof of the Borevich-Shafarevich conjecture without using resolution of singularities.

Both Igusa and Denef reduce the rationality of the Poincare series to that of $p$-adic integrals. In Igusa's theorem the integral is over $\Z_p^m$, but in the case of Denef's theorem, the domain of the integral is a definable subset of $\Q_p^m$ which does not carry any structure of a $p$-adic algebraic or analytic set. But using Macintyre's quantifier elimination, the domain is a so-called $p$-adic semi-algebraic set (analogous to the real semi-algebraic sets), and this allows computation of the integral.

Denef's results and methods of proofs have been widely used in rationality proofs for various other Poincare series. These include Poincare series counting subgroups of a group (part of the subject of subgroup growth) or representations 
of a group (part of the subject of representation growth). See Section \ref{growth} and the book by 
Lubotzky and Segal \cite{alex-dan-book} for more on this.

Denef's results have also influenced 
the study of local height zeta functions that are used in counting rational points of bounded height in algebraic varieties over number fields (via Euler products). See Chambert-Loir's paper \cite{loir-surv} and the survey \cite{adeles-surv} for more on this.

It is known that suitable quantifier elimination for $p$-adic fields yields quantifier elimination for adeles. See \cite{adeles-surv}and the references there. Furthermore, rationality of $p$-adic integrals are used in \cite{zeta1} to get results for Euler products that in turn yield results for analytic properties of adelic integrals.

We now give some details on the relevant model theory.

\section{\bf Quantifier elimination for $p$-adic fields, and uniformity in $p$}\label{qe}

Let $\cL_{rings}=\{+,.,0,1\}$ denote the language of rings. Macintyre extended this language by adding predicates for 
$n$th powers for all $n$. Let $P_n(x)$ denote the formula $\exists y \ x=y^n$. The Macintyre language is $\cL_{Mac}=\cL_{rings} \cup \{P_n(x): n\geq 2\}$.

The Macintyre language has been of fundamental importance in $p$-adic model theory since its introduction. A celebrated theorem of Macintyre states the the theory of the field of $p$-adic numbers has elimination of quantifiers in $\cL_{Mac}$. 
The elimination of quantifier is in fact proved for the theory of $p$-adically closed fields and $\Q_p$ is $p$-adically closed. 

We recall that a field is $p$-adically closed (or formally $p$-adic in the terminology of Ax-Kochen \cite{AK2}) if it 
is a model of the following axioms. These axioms can be expressed in the language $\cL_{rings}$ augmented by a predicate 
for the valuation ring (called the language of valued fields) or in the language $\cL_{rings}$ since the valuation ring is 
existentially 
$\cL_{rings}$-definable (both with parameter $p$ and without any parameter). In fact, if we take $p$ as a parameter, then an easy application of Hensel's lemma shows that $\Z_p$ is quantifier-free 
$\cL_{Mac}$-definable in $\Q_p$ using the predicate $P_2(x)$ (see \cite[Lemma 2.1]{Denefrationality}). If we do not 
take $p$ as a  parameter, then it is proved in \cite{CDLM} that there is a parameter-free existential definition of $\Z_p$ in $\Q_p$, and more generally of the valuation ring $\cO_K$ in $K$, for any given finite extension $K$ of $\Q_p$.

\

{\it Axioms for $p$-adically closed fields}

\

($i$). Sentences stating the the field has characteristic zero. 

\

($ii$). Sentences expressing the property of being a Henselian valued field $K$ as follows. Let
$f(x)=a_0x^n+\dots+a_n$ be a polynomial over the valuation ring $\cO_K$ such that 
there is some $\alpha \in \cO_K$ and $r\geq 0$ such that $v(f(\alpha))\geq 2r+1$ and $v(f'(\alpha))<r+1$. 
Then there is a unique $\bar{\alpha}\in \cO_K$ such that $f(\alpha)=0$ and $v(\alpha - \bar{\alpha}) \geq r+1$.   

\

($iii$). Sentences expressing that the value group is a $Z$-group. These sentences express that the group is abelian and has a minimal positive element $1$, and for all $n$ the index of $nG$ in $G$ is $n$ (we have written the group additively). This is equivalent to the condition that the group is elementarily equivalent to $(\Z,+,0,1,<)$ in the language of ordered abelian groups.

\

($iv$). The sentence stating that $v(p)=1$.

\

($v$). The sentence stating that the residue field is $\F_p$.

\

Ax and Kochen \cite{AK2} proved that these axioms completely characterize $p$-adically closed fields. In other words the 
theory of $p$-adically closed fields is complete in the language of rings or the language of valued fields, i.e. 
a sentence in the language holds in $\Q_p$ if and only if it holds in every $p$-adically closed field. 

Now we state Macintyre's theorem.

\begin{thm}[Macintyre \cite{Macintyre1}]\label{mac-thm} The theory of $p$-adically closed fields admits elimination of quantifiers in $\cL_{Mac}$.\end{thm}

A subset of $\Q_p^m$ that is quantifier-free definable in the Macintyre language is called a $p$-adic semi-algebraic set. Note that such a set is a Boolean combination of sets of the form 
$$\{x\in \Q_p^m: \exists y\in \Q_p: f(x)=y^n\},$$
where $f\in \Q_p[x_1,\dots,x_m]$.

Another formulation of Macintyre's theorem is that if $S\subseteq \Q_p^{m+q}$ is semi-algebraic, then 
its projection on the $m$-coordinates
$$\{x\in \Q_p^m: \exists y\in \Q_p^q, (x,y)\in S\},$$
is also semi-algebraic.

Now the question arises as to whether one can obtain a language $\cL$ in which the $\Q_p$, for all $p$ or all but finitely many $p$, have a quantifier elimination that is uniform in $p$. More generally, one can ask for a uniform quantifier elimination for 
the non-Archimedean completions $K_v$ of a number field $K$, 
where $v$ runs over non-Archimedean absolute values of $K$, for all $v$ or almost all $v$. 

If $\cL$ is many-sorted we require that $\psi(x)$ be quantifier-free of the field sort (so it could have quantifiers 
ranging over the other sorts).

We remark that a family $\{K_v: v\in S\}$, where $S$ is a subset of the set of normalized absolute values, is said to have uniform quantifier elimination if for any $\cL$-formula $\varphi(x)$, where $x$ is a tuple of variables, there exists a quantifier-free $\cL$-formula $\psi(x)$ such that for all $v\in S$ 
$$K_v \models \forall x (\varphi(x) \Leftrightarrow \psi(x)).$$

Belair \cite{Belair} and Macintyre \cite{Macintyre2} obtained uniform quantifier eliminations for $\Q_p$ for all $p$. Pas \cite{pas} obtained a quantifier elimination that is uniform for almost all $p$. Macintyre and Pas' results were used to get uniform rationality for $p$-adic integrals).

In these works, the uniform quantifier elimination for all but finitely many $\Q_p$ is deduced from a 
quantifier elimination for the theory of Henselian valued fields of characteristic zero and residue field of characteristic zero. In the Macintyre and Pas cases, this quantifier elimination is in a many-sorted language for the field sort relative to the other sorts.

\section{\bf Denef's rationality theorem and definability of the Serre series}\label{serre}

Denef's remarkable theorem on $p$-adic integrals is the following.

\begin{thm}[Denef \cite{Denefrationality}]\label{denef-thm2} Let $X$ be a definable subset of $\Q_p^m$ in the language of rings. Suppose that $X$ is contained in a compact subset of $\Q_p^m$. Let $g\in \Q_p[x]$, where $x=(x_1,\dots,x_m)$. Then
$$\int_X |g(x)|^s dx$$
is a rational function of $p^{-s}$.\end{thm}

Denef adds a predicate $|$ to $\cL_{rings}$ interpreted as $x |y \Leftrightarrow v(x) \leq v(y)$ and works in this language, 
but the relation $|$ is $\cL_{rings}$-definable using the parameter-free $\cL_{rings}$-definability of the valuation ring 
$\Z_p$ in $\Q_p$ in\cite{CDLM} (same for a finite extension of $\Q_p$). 

We note that Denef's theorem remains valid
for a finite extension $K$ of $\Q_p$, with the same proof. See \cite{denef-surv}.

In \cite{Denefrationality} more general results are proved where $g(x)$ is replaced by a definable function from $\Q_p^m$ into $\Q_p$ (i.e. a function whose graph is definable).

The reason Theorem \ref{denef-thm2} implies Serre's conjecture is that the Serre series $\tilde{P}(T)$ can be 
written as an integral as in Theorem \ref{denef-thm2} where the domain of integration is definable. 

Denef's proof of this definability result in \cite{Denefrationality} goes as follows.

Let $f_1(x),\dots,f_r(x) \in \Q_p[x]$ where $x=(x_1,\dots,x_m)$. Let $s\in \R_{>0}$. Let
$$I(s)=\int_D |w|^s dx dw,$$
where 
$$D=\{(x,w)\in \Z_p^m \times \Z_p: \exists y \in \Z_p^m, x\equiv y \ \mathrm{ mod} \ w, f_i(y)=0, \ \mathrm{for} \ i=1,\dots,r\}.$$
Then 
$$I(s)=\sum_{n\geq 0} \int_{(x,w)\in D, v(w)=n} p^{-ns} dxdw$$
$$=\sum_{n\geq 0} p^{-ns} \int_{(x,p^n)\in D, v(w)=n} dx dw$$
$$=\sum_{n\geq 0} p^{-ns} (\int_{(x,p^n)\in D} dx)(\int_{v(w)=n} dw)$$
$$=(1-p^{-1}) \sum_{n\geq 0} p^{-ns} \tilde{N}_n p^{-nm}(p^{-n}-p^{-(n+1)})$$ 
$$=(1-p^{-1})\sum_{n\geq 0} \tilde{N}_n(p^{-s}p^{-m-1})^n$$
$$=(1-p^{-1})\tilde{P}(p^{-m-1-s}).$$

This technique is powerful and can be used in a variety of contexts. A variant is done for the conjugacy class zeta function of an algebraic group and is stated in Section \ref{conj}. 

\begin{prob} Use these methods to prove definability for other generating functions in algebra and number theory.
\end{prob}

\section{\bf Subgroup growth zeta functions of groups}\label{growth}
Denef's theorem \ref{denef-thm2} led to several results on counting subgroups of a group. This started with work of 
Grunewald-Segal-Smith in \cite{GSS} who initiated the subject of subgroup growth which studies
growth of subgroups of a group $G$ of given index. 

Let $G$ be a group. One assumes that $G$ is finitely generated, so that for any $n$, it has only finitely 
many subgroups of index $n$. Let $a_n(G)$ denote this number. Grunewald-Segal-Smith 
defined the subgroup growth zeta function of $G$ as 
$$\zeta_G(s)=\sum_{n\geq 1} a_n(G) n^{-s}=\sum_{H \leq G} |G:H|^{-s}.$$
This is a non-commutative generalization of he Dedekind zeta function of a number field.

\begin{ex} Let $G=\Z^2$, the free abelian group of rank two. Then $$\zeta_{\Z^2}(s)=\zeta(s)\zeta(s-1).$$ 
This implies that 
$$a_1(G)+\dots+a_N(G) \sim (\pi^2/12)N^2$$
as $N\rightarrow \infty$. See \cite{alex-dan-book}.

\end{ex}

The function $\zeta_G(s)$ defines an analytic function of some half-plane precisely when the coefficients $a_n(G)$ 
are bounded by a polynomial. An interesting theorem of Lubotzky-Mann-Segal \cite{LMS} gave a characterization of finitely generated residually finite groups with such polynomial subgroup growth: they are precisely groups which have a subgroup 
of finite index that is soluble of finite rank.

Grunawald-Segal-Smith \cite{GSS} studied $\zeta_G(s)$ for $G$ finitely generated nilpotent. In this case one has an Euler product factorization
$$\zeta_G(s)=\prod_p \zeta_{G,p}(s)$$
where the local factors are defined by
$$\zeta_{G,p}=\sum_{n\geq 0} a_{p^n}(G)p^{-ns},$$
which is a zeta function counting subgroups of $p$-power index.

Using Denef's theorem \ref{denef-thm2}, they proved in \cite{GSS} that each local factor $\zeta_{G,p}(s)$ is a rational function of $p^{-s}$ by writing it as a $p$-adic integral over a definable subset.

The precise asymptotic growth of $a_n(G)$ for finitely generated nilpotent $G$, 
was given by du Sautoy and Grunewald in \cite{dsG}. They used Denef's theorem and some variants and other tools to prove the following remarkable result.

\begin{thm}[du Sautoy-Grunewald \cite{dsG}] \label{dsG} The abscissa of convergence of $\zeta_G(s)$ is a rational number $\alpha(G)$ and $\zeta_G(s)$ has meromorphic continuation to the half-plane $Re(s)>\alpha(G)-\delta$ for some $\delta>0$. It follows that 
$$a_1(G)+\dots+a_N(G) \sim cN^{\alpha(G)} (\mathrm{log} N)^{b(G)}$$
where $c\in \R$ and $b(G)\in \N$.\end{thm}

Their result is proved more generally for Euler products of cone integrals which are examples of Denef's integrals over definable sets. The question arises as to whether this holds for all definable integrals. This turned to be true and was proved in \cite{zeta1}. See Section \ref{euler} for more on this.

\section{\bf Conjugacy class zeta functions of algebraic groups\\ over local fields}\label{conj}

Let $\mathcal{O}$ denote a complete discrete valuation ring with maximal ideal $\cM$ and finite residue field $k$ of cardinality $q$. 
Let $G$ be a Chevalley group with an embedding into $\GL_d$. Consider the
congruence subgroups defined by 
$$G^m(\mathcal{O})=\mathrm{Ker}(G(\cO)\to \mathrm{GL}_d(\cO/\cM^m))$$
of $G(\cO)$. 

For each $m \in \N$, let $c_m$ denote the the number
of conjugacy classes in the congruence quotients 
$$G(\cO,m)\cong G(\cO)/G^m(\cO)\cong G(\cO/\cM^m).$$

du Sautoy \cite{ds-conj} defined the conjugacy class zeta function $\zeta^{\cc}_{GL_n(\Z_p)}(s)$ 
for $GL_n(\Z_p)$ and proved its rationality using Denef's Theorem \ref{denef-thm2}. It follows that the number of conjugacy classes in $GL_n(\Z_p/p^m\Z_p)$ for $m\geq 1$, and fixed $p$, satisfy a linear recurrence relation. He asked to what extent this relation depends on the prime $p$. This was answered in joint work with Mark Berman, Uri Onn, and Pirita Paajanen  \cite{BDOP} as follows. The answer is that the rationality and recurrence relations depend only on the residue field, for 
large $p$.

Extending du Sautoy's definition, we define the conjugacy class zeta function for any Chevalley group $G$ as 
$$\zeta^{\cc}_{G(\cO}(s)=\sum_{m=0}^\infty c_m q^{-ms}.$$

\begin{thm} [Berman-Derakhshan-Onn-Paajanen \cite{BDOP}]\label{BDOP} Let $G$ be a Chevalley group. Let $d=\mathrm{dim}(G)$. 
\begin{enumerate}
\item For any complete discrete valuation ring $\cO$ with residue field $k$ of cardinality $q$,
$$\zeta_{G(\cO)}^{\cc}(s)=1+\frac{q^{-\mathrm{dim}(G)}|G(k)|(Z_{G(\cO)}(s-d)-1)}{1-q^{s-d}},$$
where $$Z_{G(\cO)}(s)=\int_{G(\cO)\times G(\cO)}
||\{(\mathbf{xy-yx})_{ij}\suchthat  1\leq i,j\leq d\} ||^sd\nu$$
and $\nu$ is the normalized Haar measure on $G(\cO)\times G(\cO)$. Here
$\mathbf{x}=(x_{ij})$ and $\mathbf{y}=(y_{ij})$ are $d\times d$
matrices of indeterminates. The norm in the integrand is $||\{z_i:
i\in I\}||:=\max_{i\in I}\{|z_i| \}$.
\item There exist a constant $N\in \mathbb{N}$ and formulas 
$\psi_1(x),\dots,\psi_n(x)$ in the language of rings in an $r$-tuple of variables $x$, and rational
functions in two variables 
$$R_1(X,Y),\dots,R_n(X,Y)$$ over $\Z$ such that for all
complete discrete valuation rings $\cO$ with finite residue field
$k$ of characteristic $p$ and cardinality $q$, where
$p>N$, the conjugacy class zeta function depends only on $q$ and
can be written as
$$\zeta^{\cc}_{G(\cO}(s)=card(\psi_1(k))R_1(q,q^{-s})+\dots+card(\psi_n(k)) R_n(q,q^{-s}),$$
where $\psi_j(k)$ denotes the set defined by $\psi_i$ in $k^r$. 
\item If $\cO$ and $\cO'$ are complete discrete valuation rings with the same finite
residue field of characteristic larger than $N$, then $G(\cO)$
and $G(\cO')$ have the same conjugacy zeta function, and for all $m\geq 1$, the congruence quotients $G(\cO,m)$ and $G(\cO',m)$ have the same number of conjugacy classes.
\end{enumerate}
\end{thm}

\begin{note} \noindent
\begin{enumerate}
\item Theorem \ref{BDOP} combined with Denef's theorem \ref{denef-thm2} prove the rationality of the conjugacy class zeta function for all Chevalley groups $G$ and all $\cO$ as above.  
\item If the residue characteristic of $\cO$ is large enough then the rationality and the recurrence relations are the same for all $G(\cO)$ such that the valuation rings $\cO$ have the same residue field. This includes the positive characteristic case of $G(\F_q[[t]])$ as well. So for example $G(\Z_p)$, $G(\cO)$, for $\cO$ the valuation ring of a totally ramified extension of $\Q_p$ and $G(\F_p[[t]])$ have the same conjugacy class zeta function, for large $p$.\end{enumerate}\end{note}

These are model-theoretic transfer principles for the conjugacy class zeta functions.

\

The proof of Theorem \ref{BDOP} is based on uniform cell decomposition which is a partition of the domain of integration into cells uniformly in $p$ on each of which the integral can be evaluated. The idea of such a cell decomposition was due to Paul Cohen in his proof of decidability of $p$-adic fields in \cite{Cohen}, but that was only in dimension 1 and a general cell decomposition theorem was proved by Denef who used it to give a second proof of his rationality result in \cite{Denefrationality}. 

Uniform cell decomposition was proved by Pas in \cite{pas}. This was used in \cite{BDOP} to give a uniform in $p$ version of Denef's formulas for the integrals showing that the integrals depend only the residue field. This also follows from motivic integration results of Cluckers-Loeser \cite{CL2} and Denef-Loeser \cite{DL} but an new and self-contained treatment was given in \cite{BDOP} without introducing motivic integrals and more of an algebraic nature.

In \ref{BDOP} the rationality of the conjugacy class zeta function is proved by first proving definability of the conjugacy class zeta function, similar to Denef's proof of definability of the Serre series and du Sautoy's proof in \cite{ds-conj}. However, this definability is in the form of an integral with definable domain but with respect to the Haar measure on the 
group $G(\cO)$. To apply Denef and Pas' work one needs to work with integrals with respect to the additive Haar measure on the local field. This requires writing integrals with respect to the Haar measure on $G(\cO)$ as integrals with respect to the additive Haar measure on the local field. A general result of this kind is proved in \cite{BDOP} that should be of independent interest as follows.

\begin{thm}[Berman-Derakhshan-Onn-Paajanen\label{change}] Let $G$ be a Chevalley group defined over $\Z$. There exist (explicitly computed)
rational maps $\iota: \mathbb{A}^{\dim G} \to G$ and $\alpha:
\mathbb{A}^{\dim G} \to \mathbb{A}$ such that the following holds.
Let $K$ be a non-Archimedean local field  with valuation ring
$\cO$. Let $dx$ denote the additive Haar measure on $K^{\dim G}$
normalized on $\cO^{\dim G}$. Then the measure $\mu$ given by
$$\int_{G(K)} f(g) d\mu:=\int_{K^{\dim G}} |\alpha(x)| f(\iota(x)) dx,$$
as $f$ runs through all complex valued Borel functions on $G(K)$,
is a left and right Haar measure on $G(K)$ normalized on
$G(\cO)$.
\end{thm}

\

\section{\bf Definability of the conjugacy class zeta function: \\ Proof of Theorem \ref{BDOP}(1)}\label{def-conj}
We give the proof of definability of the conjugacy class zeta function that is joint work with Berman-Onn-Paajanen in \cite{BDOP}. 

As in Section \ref{conj}, $\cO$ denotes a complete discrete valuation ring with finite residue field $k$ of characteristic $p$, and $c_m=c_m(G,\cO)$ denotes the number of conjugacy classes in the congruence quotient $G(\cO, m)$. We have  
\begin{equation}
\zeta^{\cc}_{G(\cO)}(s)=\sum_{m=0}^{\infty} c_m~ q^{-m s}\ \ \ \
(s\in\C).
\end{equation}
We put $k_G(q)=|G(k)|q^{-d}$, $q=|k|$, and $d=\dim G$.

\begin{lemma}\cite{BDOP}\label{lemma-index-dimension-formula}
Let $X$ be a smooth scheme defined over $\cO$. Then
for all $m\geq 1$
\[
|X(\cO/\cM^m)|=|X(k)|q^{(m-1)\dim X}.
\]
In particular the formula holds for algebraic groups defined over $\cO$.

\end{lemma}

Now we give a proof of Theorem \ref{BDOP}(1).

Let $w: G(\cO)\times G(\cO)\ra \N$ be defined by 
$$w(x,y)=\max\{m\in\N: x^{-1}y^{-1}xy\in G^m(\cO)\}$$ if
$x^{-1}y^{-1}xy\neq 1$, and $\infty$ otherwise. It is
easily checked that $$w(x,y)=\min_{1\leq i,j\leq
n}\{v(xy-yx)_{ij}\}.$$
We call it the depth of a commutator in the congruence filtration

We need the following classical result for a finite group $H$: let $cc(H)$ denote the number of conjugacy classes in $H$,
then (cf. \cite{ds-conj})
\begin{equation}\label{nonburnside}
cc(H)=|H|^{-1}\sum_{x\in H} |C_H(x)|= |H|^{-1}|\{(x,y)\in H\times H:xy=yx  \}|.
\end{equation}

Applying \eqref{nonburnside} to $G(\cO,m)$ we see that for
all $m$, 
$$c_m=|G(\cO)/G^m(\cO)|^{-1}e_m,$$ 
where

\[e_m:=|\{({x},{y})\in G( \cO,m)\times G(\cO,m):
{x}^{-1}y^{-1}xy=1\}|.\]

Let

\[
\begin{split}
&W_m =\{({x},{y})\in G(\cO)\times G(\cO): w({x},{y})\geq m \},\\
&Z'_{G(\cO)}(s)=\sum_{m=0}^\infty q^{-ms}\nu(W_m).
\end{split}
\]
Then
\begin{align*}Z_{G(\cO)}(s)& = \sum_{m=0}^\infty q^{-ms}\nu(W_m \smallsetminus W_{m+1})=
Z'_{G(\cO)}(s)(1-q^s)+q^s.
\end{align*}
It is easily seen that 
$$\nu(W_m)=\nu(G^m(\cO)\times G^m(\cO))e_m.$$ 
Applying Lemma~\ref{lemma-index-dimension-formula} yields

\begin{align*}
Z'_{G(\cO)}(s)&=\sum_{m=0}^{\infty} \nu (G^m(\cO)\times
G^m(\cO)) e_m q^{-ms}\\
&= \sum_{m=0}^{\infty} |G(\cO)/G^m(\cO)|^{-1}c_mq^{-ms}\\
&=1+k_G(q)^{-1}(\zeta_{G(\cO)}^{\cc}(s+d)-1),
\end{align*}

since $\nu(G^m(\cO)\times G^m(\cO))=|G(\cO)/ G^m(\cO)|^{-2}$. The proof is complete.

\section{\bf From $p$-adic integrals to global zeta functions\\ via Euler products}\label{euler}

After Denef proved his rationality Theorem \ref{denef-thm2}, Pas \cite{pas} and Macintyre \cite{Macintyre2} independently proved that there are uniformities in the shape of the rational functions as $p$ varies when $X$ and $g(x)$ are over $\Q$. 

When there such a uniformity in local $p$-adic integrals, one expects that it is explained by some 
global information of number-theoretic nature. By global we mean relating to a global field that is a finite extension of 
$\Q$ (called a number field) or a field of rational functions of a curve over a finite field $\F_q((t))$.

In 1995 Kontsevich gave a seminal talk at Orsay where he introduced the idea of motivic integration \cite{maxim-talk}. He used this to prove that birationally equivalent Calabi-Yau manifolds have the same Hodge numbers. Batyrev had proved they have the same Betti numbers using $p$-adic integration and the Weil conjectures. Kontsevich avoided this by replacing $p$-adic integration by motivic integration that comes from a measure that is not $\R$-valued but valued in a Grothendieck ring of algebraic varieties, and did gain the extra important geometric information. See Loeser's Seattle lectures \cite{francois-surv} for more details and an introduction to motivic integration.

In a series of works, Denef and Loeser developed a systematic theory of motivic integration, see \cite{DL} and the survey 
\cite{francois-surv}. This theory has had various applications to algebraic geometry, number theory, and representation theory. On the other hand it gave a beautiful geometric explanation of the uniformity in $p$ of the rational functions by showing that there is a "motivic integral" that relates to a "Chow motive" in the sense of algebraic geometry that specializes via a certain procedure to the $p$-adic integrals for almost all $p$, see \cite{DL}. 

A new approach to motivic integration was given by Cluckers-Loeser, see for example \cite{CL2}, and Hrushovski-Kazhdan \cite{HK}. Each work has had numerous applications to geometry and arithmetic.

A different approach to use the uniformity in the $p$-adic integrals and get global information from it was introduced in \cite{zeta1} via the idea of taking an Euler product. Several Dirichlet series in number theory and geometry are Euler products, for example the Riemann zeta function
$$\sum_{n\geq 1} n^{-s}=\prod_p (1-p^{-s})^{-1}$$
and various other zeta functions and $L$-functions, e.g. $L$-functions of algebraic varieties and $L$-functions of Galois representations and modular forms. Furthermore, integrals over the ring of adeles $\A_K$ of a number field 
are naturally Euler products. These $L$-functions play a fundamental and central role in modern number theory and representation theory. See \cite{adeles-surv} for more on these.

The question then arises as to whether an Euler product over all primes $p$ of the $p$-adic integrals (suitably normalized)
over definable sets (the kind considered by Denef in his rationality theorem) would have good analytic properties similar to zeta and $L$-function, for example meromorphic continuation and information on its poles. 

Much light on this problem is shed by the work of du Sautoy and Grunewald in \cite{dsG} in the proof of Theorem 
\ref{dsG} on the subgroup growth zeta function of finitely generated nilpotent groups. In fact, in \cite{dsG} they prove 
the conclusions of Theorem \ref{dsG} for an Euler product of "cone integrals". These cone integrals are a special case of the definable integrals of Denef. 

In \cite{zeta1}, Theorem \ref{dsG} of du Sautoy-Grunewald was generalized to all "definable integrals" as follows.

\begin{thm}[Derakhshan {\cite{zeta1}}]\label{Thm-zeta} Let $Z(s,p)$ be as above. 
Let $a_{p,0}$ be the constant coefficient of $Z(s,p)$ when expanded as a power series in $q^{-s}$. Then the Euler product over all primes $p$ 
$$\prod_p a_{p,0}^{-1} Z(s,p)$$
has rational abscissa of convergence $\alpha$ and meromorphic continuation to the half-plane $\{s: Re(s)>\alpha -\delta\}$ for some $\delta>0$.
The continued function is holomorphic on the line $Re(s)=\alpha$ except for a pole at $s=\alpha$.\end{thm}
Tauberian theorems from analytic number theory then imply the following asymptotic formulas.
\begin{cor}[Derakhshan {\cite{zeta1}}] \label{cor-zeta} 
Suppose the Euler product $\prod_p a_{p,0}^{-1} Z(s,p)$ can be written as a Dirichlet series $\sum_{n\geq 1} a_n n^{-s}$. Then for some real numbers $c,c'\in \R$, 
$$a_1+a_2+\dots+a_N \sim c N^{\alpha}(log N)^{w-1}$$ 
$$a_1+a_22^{-\alpha}+\dots + a_N N^{-\alpha} \sim c'(log N)^w$$
as $N \rightarrow \infty$, were $w$ is the order of the pole of $Z(s,p)$ at $\alpha$.\end{cor}

du Sautoy has given examples of subgroup zeta functions that are of the form of the Euler products 
in Theorem \ref{Thm-zeta} but do not admit meromorphic continuation to the entire complex plane. See \cite{ds-nat-bnd}, 
\cite{ds-book}.

We remark that in the proof of Theorem \ref{Thm-zeta}, we found it necessary to use the uniform formulas of 
Macintyre for the $p$-adic integrals that are proved in \cite{Macintyre2}. In this work, Macintyre introduces a many-sorted language for $\Q_p$ for all $p$ uniformly in $p$, and proves uniform quantifier elimination and cell decomposition. In most uniform treatments 
of $\Q_p$. The fact that this is proved for all $p$ and the nature of the sorts in the Macintyre language for each $\Q_p$ is crucial in order to prove Theorem \ref{Thm-zeta}.

Since most integrals over the ring of adeles $\A_K$ of a number field $K$ are such Euler products, 
Theorem \ref{Thm-zeta} and Corollary \ref{cor-zeta} have applications to adelic integrals and number theory. See \cite{adeles-surv} for such results and problems and relations to various $L$-functions in number theory and the Langlands conjectures.

A challenging problem stated in \cite{adeles-surv} asks for adding "Archimedean factors" to these products, to obtain a function similar to a zeta or $L$-function in adelic form.

\section{\bf Global conjugacy class zeta function}\label{global-conj}

Uri Onn defined the global conjugacy class zeta function of an algebraic group. It counts conjugacy classes in algebraic groups over a number field, and can be understood via its local factors using \ref{Thm-zeta} and the definability of the local conjugacy class zeta functions in \cite{BDOP} discussed in Section \ref{def-conj}.

For simplicity we consider the case $K=\Q$.
For the definition and basic properties of strong approximation in algebraic groups we refer to the book by Platonov and Rapinchuk \cite[Chapter 7.1]{Platonov-R-book}. For example, $SL_n$ has strong approximation but $GL_n$ does not, for all $n$. 

Let $G \le \GL_n$ be a
$\Z$-defined algebraic subgroup. Let $c_m$ denote 
the number of conjugacy classes in $G(\Z/m\Z)$.
The global conjugacy zeta function of $G(\Z)$ is

$$\zeta_{G(\Z)}^{\mathrm{cc}}(s)=\sum_{m\geq 1}
c_{m} m^{-s}$$

\begin{thm}\label{glob-conj} \cite{zeta1} Assume that $G$ has strong approximation. 
Then the global conjugacy zeta function $\zeta_{G(\Z)}^{\mathrm{cc}}(s)$ has rational abscissa of convergence $\alpha$ and meromorphic continuation to the half-plane
$\{s: Re(s)>\alpha-\delta\}$ for some $\delta >0$. The continued function is holomorphic on the line $Re(s)=\alpha$ except for a pole at $s=\alpha$ of some order $w$. 
There exists $c\in \R$ such that 
$$c_1+\dots+c_N \sim cN^{\alpha} (log N)^{w-1}$$
as $N\rightarrow \infty$.\end{thm}

To see that Theorem \ref{glob-conj} follows from Theorem \ref{Thm-zeta} note that by strong approximation, for any $m=p_1^{r_1}\dots p_k^{r_k}$, we have an isomorphism

$$G(\Z/m\Z) \cong G(\Z/p_1^{r_1}\Z) \times \dots
\times G(\Z/p_k^{r_k}\Z).$$

For the irreducible representations (which are the same as conjugacy classes) we have identifications

$$\mathrm{Irr}(G(\Z/m\Z)) =
\mathrm{Irr}(G(\Z/p_1^{r_1}\Z))\times \dots
\times \mathrm{Irr}(G(\Z/p_k^{r_k})).$$

Therefore, the numbers of conjugacy classes are multiplicative, and we deduce an Euler factorization 
$$\zeta_{G(\Z)}^{\mathrm{cc}}(s)=\prod_{p}\zeta_{G(\Z_p)}^{\mathrm{cc}}(s).$$

\begin{prob}\label{conj-zeta} Complete the Dirichlet series $\zeta_{G(\Z)}^{\mathrm{cc}}(s)$ by adding "Archimedean factors" and write it as an integral of a suitable function over a definable subset of $\A_K^m$ for some $m\geq 1$. 
Understand the analytic properties (e.g. meromorphic continuation and functional equation) of the completed zeta function.
\end{prob}
For more on adelic integration and related problems see \cite{adeles-surv}.

\section{\bf Zeta functions arising from definable equivalence relations\\ and $p$-adic elimination of imaginaries}

Motivated by counting irreducible complex representations of a finitely generated nilpotent group up to equivalence and iso-twisting (see \cite{HMR}), 
Hrushovski, Martin, and Rideau \cite{HMR} proved an extension of Denef's Theorem \ref{denef-thm2} to generating functions that count the number of  classes of a parametrized family of definable equivalence relation. We state this result.

We consider $\Z$ in the language of ordered Abelian groups and $\Q_p$ in the language of rings. We consider $\Z$ as the value group of $\Q_p$ and add a sort for it together with a symbol $v$ for the valuation from $\Q_p$ into $\Z\cup \{\infty\}$.  By a definable family $R=(R_l)_{i\in \Z^r}$ of subsets of $\Q_p^N$ we mean a definable subset 
$R\subseteq \Q_p^N \times \Z^r$, where $R_l$ denotes the fibre above $l$ of the projection map $R\rightarrow \Z^r$. 
By a definable family of equivalence relations we mean an equivalence relation $E=(E_l)_{l\in \Z^r}$ on $R$ such that 
for any $x,y$ if $xEy$ then there is some $l\in \Z^r$ such that $x,y \in R_l$. This induces a definable equivalence relation 
$E_l$ on $R_l$ for each $l$. 

We may regard $(E_l)_{l\in \Z^r}$ as a definable family of subsets of $\Q_p^{2N}$. Since $\N_0$ (the non-negative integers) is definable in $\Z$, we may consider the definable families $R=(R_l)_{i\in \N_0^r}$. 

Let $t=(t_1,\dots,t_r)$. For $l\in \N^r$, we let $t^l$ denote $\prod_{1\leq i \leq r} t_i^{l_i}$. 

\begin{thm}[Hrushovski-Martin-Rideau \cite{HMR}]\label{HMR} Let $R=(R_l)_{l\in \N_0^r}$ be a definable family of subsets of $\Q_p^N$. Let $E=(E_l)_{l\in \N_0^r}$ be a definable family of equivalence relations on $(R_l)_{l\in \N_0^r}$. Suppose that for all $l\in \N_0^r$ the set of equivalence classes $R_l/E_l$ is finite of cardinality $a_l$. Then the Poincare series $\sum_{l\in \N_0^r}a_l t^l \in \Q[[t_1,\dots,t_r]]$ is a rational function.\end{thm}
We remark that the series $\sum_{l\in \N_0^r}a_l t^l \in \Q[[t_1,\dots,t_r]]$ is said to be a rational function if it is 
equal to a rational function in $t_1,\dots,t_r$ with coefficients from $\Q$.

 Cluckers \cite[Appendix]{HMR} gave a different proof of Theorem \ref{HMR} in the more general context with  subanalytic language of $\Q_p$. 

In order to prove Theorem \ref{HMR}, Hrushovski-Martin-Rideau prove that the theory of $\Q_p$ admits elimination of 
imaginaries in the geometric language. This language has a sort for $\Q_p$ with the language of rings and sorts for the
quotients $S_n=GL_n(\Q_p)/GL_n(\Z_p)$ with the language of groups and symbols for the projection maps $GL_n(\Q_p) \rightarrow S_n$. Note that $S_1$ can be identified with the value group and the corresponding projection map with the 
valuation map. Moreover, they prove that the elimination of imaginaries in uniform in $p$. They also prove their results for finite extensions of $\Q_p$. See \cite{HMR} for details.

Theorem \ref{HMR} applies naturally to problems in representation growth of groups. This is a recently
developing subject, analogous to subgroup growth, which studies the numbers of irreducible complex representations in each 
degree using representation zeta functions.

Let $G$ be a finitely generated group having finitely many irreducible representations of any fixed dimension, up to equivalence (such groups are called representation rigid). Let $r_n(G)$ denote this number. Let
$$\zeta_{G}^{Rep}(s)=\sum_{n\geq 1} r_n(G)n^{-s}=\sum_{\rho\in \mathrm{Irr}(G)}(dim(\rho))^{-s},$$
where $\mathrm{Irr}(G)$ denotes the set of finite dimensional, complex, and irreducible representations of $G$.

If $G$ is a finitely generated nilpotent group, then $G$ is not representation rigid, but if we take iso-twist classes of representations in each degree which are classes of the equivalence relation between representation when one twists by a 1-dimensional representation, then we get a finite set and thus one can define a representation zeta function. 
Hrushovski-Martin-Rideau apply Theorem \ref{HMR} to prove rationality of these zeta functions. See \cite{HMR} for details. 

If the sequence $r_n(G)$ grows polynomially, then $\zeta_{G}^{Rep}(s)$ converges in some half-plane $\{s: Re(s)>\alpha\}$ 
and we may speak of its abscissa of convergence $\alpha_{G}$.

Avni \cite{Nir} proved the following

\begin{thm}[Avni \cite{Nir}] Let $\Gamma$ be an arithmetic lattice in characteristic zero (i.e. a lattice in a higher semi-simple higher rank Lie group) with the congruence subgroup property. Then $\alpha_{\Gamma}$ is a rational number.\end{thm}
In the proof, $\zeta_{\Gamma}^{Rep}(s)$ is written as an Euler product of certain $p$-adic integrals over definable sets which are more general than the integrals in Denef Theorem \ref{denef-thm2}, which are also rational functions. Then an estimate is given for their abscissa of convergence.

\bibliographystyle{acm}
\bibliography{bibadeles}

\end{document}